\documentclass[11pt,a4paper]{article}
\title{\bf Fokker-Planck type equations with Sobolev diffusion coefficients
and BV drift coefficients}
\author{Dejun Luo\footnote{Email: luodj@amss.ac.cn. This work is supported by the
Key Laboratory of Random Complex Structures and Data Science, Academy of Mathematics
and Systems Science, Chinese Academy of Sciences (No. 2008DP173182).}
\vspace{3mm}\\
{\footnotesize Institute of Applied Mathematics, Academy of Mathematics and Systems Science,}\\
{\footnotesize Chinese Academy of Sciences, Beijing 100190, China} }
\date{}
\usepackage{amssymb,amsmath,amsfonts,amsthm,array,color}

\def\R{\mathbb{R}}

\def\D{\mathcal{D}}
\def\K{\mathcal{K}}

\newcommand{\ra}{\rightarrow}
\newcommand{\lra}{\longrightarrow}
\newcommand{\da}{\downarrow}
\newcommand{\ua}{\uparrow}
\newcommand{\ee}{\varepsilon}
\def\BV{\textup{BV}}
\def\div{\textup{div}}
\def\d{\textup{d}}
\def\supp{\textup{supp}}

\def\fin{\hfill$\square$}
\def\<{\langle}
\def\>{\rangle}

\begin{document}

\maketitle
\makeatletter 
\renewcommand\theequation{\thesection.\arabic{equation}}
\@addtoreset{equation}{section}
\makeatother 

\newtheoremstyle{newthm}
 {3pt}
 {3pt}
 {\itshape}
 {}
 {}
 {\textbf{.}}
 {.5em}
 {\thmname{\textbf{#1}}\thmnumber{ \textbf{#2}}\thmnote{ {\textbf{(#3)}}}}

\theoremstyle{newthm}

\newtheorem{theorem}{Theorem}[section]
\newtheorem{lemma}[theorem]{Lemma}       
\newtheorem{corollary}[theorem]{Corollary}
\newtheorem{proposition}[theorem]{Proposition}
\newtheorem{remark}[theorem]{Remark}
\newtheorem{example}[theorem]{Example}
\newtheorem{definition}[theorem]{Definition}

\begin{abstract}
In this paper we give an affirmative answer to an open question
mentioned in [Le Bris and Lions, Comm. Partial Differential
Equations 33 (2008), 1272--1317], that is, we prove the
well-posedness of the Fokker-Planck type equations with Sobolev
diffusion coefficients and BV drift coefficients.
\end{abstract}



\section{Introduction}

The celebrated Di Perna-Lions theory, which was initiated in
\cite{DiPernaLions89}, claims that the well-posedness of the
transport equation implies the existence and uniqueness of the
quasi-invariant flow generated by a Sobolev vector field, provided
its divergence is bounded. See \cite{CiprianoCruzeiro05} for the
case of Gaussian measure as the reference measure. This theory has
subsequently been extended to the case of BV vector fields by
Ambrosio \cite{Ambrosio04, Ambrosio08}, to the infinite dimensional
Wiener space in \cite{AmbrosioFigalli09, FangLuo10}, and to the case
of SDE in \cite{Zhang10, FangLuoThalmaier}, inspired by the work of
Crippa and de Lellis \cite{CrippadeLellis}. In two recent papers
\cite{LeBrisLions04, LeBrisLions08}, Le Bris and Lions used again
the ideas of the Di Perna-Lions theory to study the Fokker-Planck
type equations, showing the existence and uniqueness of solutions
when the coefficients have a certain Sobolev regularity (see also
\cite{Figalli, RocknerZhang10} for some related results). In
\cite{Luo10}, the author studied the well-posedness of Fokker-Planck
type equations on the Wiener space, mainly under the exponential
integrability of the gradients and divergences of the coefficients.
For the study of Fokker-Planck equations in the Hilbert space, see
\cite{BogachevPratoRockner09} and the references therein.

The present work is motivated by \cite{LeBrisLions08}. We consider
the following Fokker-Planck equations
  \begin{equation}\label{FPE}
  \partial_t p+\partial_i(p b_i )-\partial_{ij}^2(a_{ij}p)=0,
  \end{equation}
where the coefficients $b=(b_1,\cdots,b_n)$ and $a=(a_{ij})_{1\leq
i,j\leq n}$ have only weak regularity on the spatial variables, e.g.
Sobolev or BV regularity. For the motivation of studying this kind
of equations with irregular coefficients, see \cite[Section
4]{LeBrisLions08}. The above equation \eqref{FPE} is closely related
to the stochastic differential equation
  $$\d X_t=\sigma(X_t)\,\d W_t+b(X_t)\,\d t,\quad X_0=x,$$
in which $\sigma=(\sigma_{ik})_{1\leq i\leq n,1\leq k\leq m}$ is a
matrix-valued function such that $a=\frac12 \sigma\sigma^\ast$ (here
$\sigma^\ast$ is the transpose of $\sigma$), and $W_t$ is an $m$-dimensional
standard Brownian motion. In the following we
always assume that the coefficient $a$ has such a form, as in
\cite{LeBrisLions08}. The adjoint equation of \eqref{FPE}, called
the backward Kolmogorov equation, reads as
  \begin{equation}\label{Kolmogorov}
  \partial_t p-b_i \partial_i p-a_{ij}\partial_{ij}^2 p=0.
  \end{equation}

We will mainly focus on the Fokker-Planck equation of divergence
form (a notion introduced in (5.8) of \cite{LeBrisLions08}):
  \begin{equation}\label{FPE.1}
  \partial_t p+\partial_i(pb_i)-\frac12\partial_i(\sigma_{ik}\sigma_{jk}
  \partial_j p)=0.
  \end{equation}
As pointed out in \cite[Subsection 5.1]{LeBrisLions08}, this
equation is relatively easier to deal with, compared to the
Fokker-Planck equation \eqref{FPE} and the backward Kolmogorov
equation \eqref{Kolmogorov}, since its second order term
$\partial_i(\sigma_{ik}\sigma_{jk}\partial_j p)$ has a self-adjoint
form. To state the main result of this paper, we introduce some
notations. For a locally integrable function $f:\R^n\ra\R$, we say
that $f$ has $\BV_{loc}$ regularity, denoted by $f\in
\BV_{loc}(\R^n)$, if for every $i=1,\cdots,n$, the distributional
partial derivative $\partial_i f$ of $f$ is a locally finite Radon
measure on $\R^n$. $f^-:=-(f\wedge 0)$ denotes the negative part of
the function $f$. For a $\BV_{loc}$ vector field $b:\R^n\ra\R^n$, we write
$D\cdot b=\sum_{i=1}^n \partial_i b_i$ for its divergence,
which is still a Radon measure. If $D\cdot b$ is absolutely
continuous with respect to the Lebesgue measure $\d x$, we denote by
$\div(b)$ its Radon-Nikodym density. Our main result is the following

\begin{theorem}\label{Main_thm}
Assume that $b$ and $\sigma$ satisfy:
  \begin{equation}
  \begin{array}{c}
  \displaystyle b\in \big(L^1([0,T],\BV_{loc}(\R^n))\big)^n,\quad
  \frac b{1+|x|}\in \big(L^1([0,T],L^1+L^\infty(\R^n))\big)^n,\\
  \displaystyle \div(b)\in L^1\big([0,T],L^1_{loc}(\R^n)\big),
  \quad [\div(b)]^-\in L^1\big([0,T],L^\infty(\R^n)\big);
  \end{array}
  \end{equation}
and
  \begin{align}
  \sigma\in \big(L^2([0,T],W^{1,2}_{loc}(\R^n))\big)^{n\times
  m},\quad
  \frac \sigma{1+|x|}\in \big(L^2([0,T],L^2+L^\infty(\R^n))\big)^{n\times
  m}.
  \end{align}
Then for each initial condition in $L^2\cap L^\infty$ (resp.
$L^1\cap L^\infty$), the equation \eqref{FPE.1} has a unique
solution in the space
  $$\{p\in L^\infty([0,T],L^2\cap L^\infty)\ (resp.\
  L^\infty([0,T],L^1\cap L^\infty)),\ \sigma^\ast\nabla p\in L^2([0,T],L^2)\}.$$
\end{theorem}

This theorem gives a positive answer to the open question raised in
\cite[Remark 12]{LeBrisLions08} on page 1299. It also generalizes
\cite[Theorem 1.4]{Figalli}, since the diffusion coefficient $a$ in
\cite{Figalli} is assumed to be independent on the spatial variables.
The proof of Theorem 1.1 will be provided in the next section. As in
\cite[Section 7]{LeBrisLions08}, we can adapt the above result to
the other Fokker-Planck type equations \eqref{FPE} and
\eqref{Kolmogorov}, by imposing suitable conditions on
$b^\sigma=b-\frac12 \div(\sigma\sigma^\ast)$. Finally we give a
brief discussion on the reason why we cannot further generalize it
to the case where $\sigma$ has only BV regularity, see Remark
\ref{sect-2-rem.2}.

\section{Proof of Theorem \ref{Main_thm}}

In this section, we present the proof of Theorem \ref{Main_thm}. We
consider the Fokker-Planck equation of divergence form
  \begin{equation}\label{FPE_divergence}
  \partial_t p+\partial_i(pb_i)-\frac12\partial_i(\sigma_{ik}\sigma_{jk}
  \partial_j p)=0.
  \end{equation}
First we give the mathematical meaning of the above equation. For a
given initial condition $p_0\in L^2\cap L^\infty$ (resp. $L^1\cap L^\infty$), a function
$p\in L^\infty([0,T],L^2\cap L^\infty)$ (resp.
$L^\infty([0,T],L^1\cap L^\infty)$) satisfying
$\sigma^\ast\nabla p\in L^2([0,T],L^2)$ is called a weak
solution to \eqref{FPE_divergence} if for all $\varphi\in C_c^\infty\big(
[0,T)\times\R^n\big)$, it holds
  $$\int_0^T\!\!\int_{\R^n}p\,\partial_t \varphi\,\d x\d t+\int_{\R^n}p_0\varphi(0,\cdot)\,\d x
  =-\int_0^T\!\!\int_{\R^n}p\, \<b,\nabla\varphi\>\,\d x\d t+\frac12\int_0^T\!\!\int_{\R^n}
  \<\sigma^\ast\nabla p,\sigma^\ast\nabla\varphi \>\,\d x\d t.$$

The existence of solutions to the above equation is the easier part,
see for instance the beginning of \cite[Subsection
5.4]{LeBrisLions08} for the case $p_0\in L^2\cap L^\infty$, and
\cite[Subsection 6.1.2]{LeBrisLions08} for the case $p_0\in L^1\cap
L^\infty$. Hence in the following we focus on the uniqueness part of
Theorem \ref{Main_thm} and follow the ideas in \cite[Subsection
5.3]{LeBrisLions08}. The main difference is that a single
convolution kernel is not enough to achieve our result, instead, we
need a family of kernels as in the proof of \cite[Theorem 3.5]{Ambrosio04},
see also \cite[Theorem 5.1]{Ambrosio08}.

Let
  $$\K=\bigg\{\rho\in C_c^\infty(\R^n,\R_+):\supp(\rho)\subset B(1),
  \int_{\R^n}\rho\,\d x=1\bigg\}$$
be the family of candidate convolution kernels, where $B(1)$ is the open
unit ball centered at the origin $0$. Take $\rho\in\K$ and
define $\rho_\ee=\ee^{-n}\rho(\ee^{-1}\cdot)$ for $\ee>0$. We
regularize the equation \eqref{FPE_divergence} in the spatial
variables:
  \begin{equation}\label{sect-2.1}
  \partial_t(\rho_\ee\ast p)+\rho_\ee\ast \partial_i(pb_i)
  -\frac12\rho_\ee\ast \partial_i(\sigma_{ik}\sigma_{jk}
  \partial_j p)=0.
  \end{equation}
We denote by $p_\ee=\rho_\ee\ast p$ and introduce the notation
  $$[\rho_\ee,c](f)=\rho_\ee\ast(cf)-c(\rho_\ee\ast f)$$
for a differential operator $c$. Note that $c$ can also be a real
valued function. Using this notation, we have
  \begin{align}\label{sect-2.2}
  \rho_\ee\ast\partial_i(pb_i)
  &=\rho_\ee\ast(\div(b)p)+\rho_\ee\ast(b_i\partial_ip)\cr
  &=[\rho_\ee,\div(b)](p)+[\rho_\ee,b_i\partial_i](p)
  +\partial_i(b_ip_\ee)\cr
  &=Q_{1,\ee}+Q_{2,\ee}+\partial_i(b_ip_\ee),
  \end{align}
where we have defined
  \begin{equation}\label{sect-2.3}
  Q_{1,\ee}=[\rho_\ee,\div(b)](p)\quad\mbox{and}
  \quad Q_{2,\ee}=[\rho_\ee,b_i\partial_i](p).
  \end{equation}
In fact, it is the term $Q_{2,\ee}$ that causes the trouble in the
BV situation and marks the difference between the present work and
\cite{LeBrisLions08}.

Next we have
  \begin{align*}
  \rho_\ee\ast \partial_i(\sigma_{ik}\sigma_{jk}\partial_j p)
  &=\rho_\ee\ast \big((\partial_i\sigma_{ik})\sigma_{jk}\partial_j p
  +\sigma_{ik}\partial_i(\sigma_{jk}\partial_j p)\big)\cr
  &=[\rho_\ee,\partial_i\sigma_{ik}](\sigma_{jk}\partial_j p)
  +(\partial_i\sigma_{ik})\,\rho_\ee\ast(\sigma_{jk}\partial_j p)\cr
  &\hskip12pt +[\rho_\ee,\sigma_{ik}\partial_i](\sigma_{jk}\partial_j p)
  +\sigma_{ik}\partial_i(\rho_\ee\ast(\sigma_{jk}\partial_j p)).
  \end{align*}
Define
  \begin{equation}\label{S-epsilon_T-epsilon}
  S_\ee=[\rho_\ee,\partial_i\sigma_{ik}](\sigma_{jk}\partial_j p)\quad
  \mbox{and}\quad T_\ee=[\rho_\ee,\sigma_{ik}\partial_i](\sigma_{jk}\partial_j p).
  \end{equation}
We remark that the term $S_\ee$ prevents us from extending the
results to the case where the diffusion coefficient $\sigma$ has
only BV regularity, see Remark \ref{sect-2-rem.2} for more details.
Now we have
  $$\rho_\ee\ast \partial_i(\sigma_{ik}\sigma_{jk}\partial_j p)
  =S_\ee+T_\ee+\partial_i\big(\sigma_{ik}\rho_\ee\ast(\sigma_{jk}\partial_j p)\big).$$
If we denote by
  \begin{equation}\label{commutator.1}
  R_{k,\ee}=[\rho_\ee,\sigma_{jk}\partial_j](p),\quad k=1,\cdots,m;
  \end{equation}
then
  $$\rho_\ee\ast(\sigma_{jk}\partial_j p)=R_{k,\ee}+\sigma_{jk}\partial_j p_\ee.$$
Summing up the above discussions, we arrive at the equality (5.17)
in \cite{LeBrisLions08}:
  \begin{equation}\label{sect-2.4}
  \rho_\ee\ast \partial_i(\sigma_{ik}\sigma_{jk}\partial_j p)
  =S_\ee+T_\ee+\partial_i(\sigma_{ik}R_{k,\ee})
  +\partial_i(\sigma_{ik}\sigma_{jk}\partial_j p_\ee).
  \end{equation}
Combining \eqref{sect-2.1}, \eqref{sect-2.2} and \eqref{sect-2.4},
we get an equation of $p_\ee$ similar to \eqref{FPE_divergence}, but
with some error terms on the right hand side:
  \begin{equation}\label{FPE_divergence_error}
  \partial_t p_\ee+\partial_i(p_\ee b_i)-\frac12\partial_i(\sigma_{ik}\sigma_{jk}
  \partial_j p_\ee)=-Q_{1,\ee}-Q_{2,\ee}+\frac12\big(\partial_i(\sigma_{ik}R_{k,\ee})+S_\ee+T_\ee\big).
  \end{equation}

Now we need the classical commutator estimate in the Di Perna-Lions
theory (see \cite[Lemma II.1]{DiPernaLions89} or \cite[Lemma
1]{LeBrisLions08}). We include it here for the sake of the readers'
convenience.

\begin{lemma}[Commutator estimate I]\label{commutator-estimate}
For $r,\alpha,r_1,\alpha_1\geq1$, set
$\frac1\beta=\frac1r+\frac1\alpha$ and
$\frac1{\beta_1}=\frac1{r_1}+\frac1{\alpha_1}$. Let $f\in L^{r_1}
\big([0,T],L^r_{loc}(\R^n)\big)$, $g\in
L^{\alpha_1}\big([0,T],L^\alpha_{loc}(\R^n)\big)$ and $c\in
\big(L^{\alpha_1} \big([0,T],W^{1,\alpha}_{loc}(\R^n)\big)\big)^n$.
Then as $\ee\ra0$,
  \begin{equation}\label{commutator-estimate.1}
  [\rho_\ee,c_i\partial_i](f)\ra 0\quad \mbox{in }
  L^{\beta_1}\big([0,T],L^\beta_{loc}(\R^n)\big),
  \end{equation}
and
  \begin{equation}\label{commutator-estimate.2}
  [\rho_\ee,g](f)\ra 0\quad \mbox{in }
  L^{\beta_1}\big([0,T],L^\beta_{loc}(\R^n)\big).
  \end{equation}
\end{lemma}

\begin{remark}\label{sect-2-rem.1}
{\rm We observe that under the assumptions of the above lemma,
for any fixed $\ee>0$, the commutator $[\rho_\ee,c_i\partial_i](f)$ belongs to the space
$L^{\beta_1}\big([0,T],W^{1,\beta}_{loc}(\R^n)\big)$, i.e., it has the first order
Sobolev regularity with respect to the spatial variables. This can be seen from
its expression (see (3.8) in \cite{Ambrosio04} or
\cite[Lemma 2.5]{Fang09}):
  \begin{align*}
  [\rho_\ee,c_i\partial_i](f)(x)
  &=-[(f\div(c))\ast\rho_\ee](x)+\int_{\R^n}f(y)\<c(y)-c(x),(\nabla\rho_\ee)(x-y)\>\,\d y\cr
  &=-[(f\div(c))\ast\rho_\ee](x)+[(fc)\ast(\nabla\rho_\ee)](x)
  -\<c,f\ast(\nabla\rho_\ee)\>(x).
  \end{align*}
Indeed, the three terms $(f\div(c))\ast\rho_\ee,\, (fc)\ast(\nabla\rho_\ee)$
and $f\ast(\nabla\rho_\ee)$ are smooth, therefore the assertion follows from the fact that
$c\in \big(L^{\alpha_1} \big([0,T],W^{1,\alpha}_{loc}(\R^n)\big)\big)^n$. \fin }
\end{remark}

We concentrate in the following on the $L^2$-theory, that is, we
prove the uniqueness of solutions to \eqref{FPE_divergence} in the
space
  \begin{equation*}
  X_2=\big\{p\in L^\infty\big([0,T],L^2\cap L^\infty\big);
  \sigma^\ast\nabla p\in L^2\big([0,T],L^2\big)\big\}.
  \end{equation*}
For the $L^1$-theory, the same argument as in \cite[Subsection
6.1.2]{LeBrisLions08} works.

Now we consider the error terms in \eqref{FPE_divergence_error}.
Using Lemma \ref{commutator-estimate}, we see that for the time-dependent
vector field $b\in \big(L^1([0,T], \BV_{loc}(\R^n))\big)^n$ satisfying
$\div(b)\in L^1\big([0,T],L^1_{loc}(\R^n)\big)$, the first error term
  \begin{equation}\label{Q1}
  Q_{1,\ee}=[\rho_\ee,\div(b)](p)\stackrel{\ee\ra0}\lra 0
  \quad \mbox{in } L^1\big([0,T],L^1_{loc}(\R^n)\big),
  \end{equation}
since $p\in L^\infty\big([0,T],L^\infty(\R^n)\big)$.

The estimate of
the commutator $Q_{2,\ee}$ cannot be obtained from Lemma \ref{commutator-estimate},
instead, we will rely on the work of Ambrosio (see \cite[Theorem 3.2]{Ambrosio04}
or \cite[Section 5]{Ambrosio08}). We first introduce some notations. Let
$D b_t=(\partial_j b_i(t))_{1\leq i,j\leq n}$ be the ``Jacobi'' matrix of
distributional derivatives of $b_t$, whose entries are locally finite Radon
measures since $b_t\in \BV_{loc}(\R^n)$. Denote by $|D b_t|$ its total variation.
Let
  $$D b_t=D^a b_t+D^s b_t$$
be the Lebesgue decomposition of $Db_t$ into absolutely continuous and
singular part with respect to the Lebesgue measure $\d x$. Define the matrix
valued function $M:[0,T]\times\R^n\ra\R^n\otimes\R^n$ as the Radon-Nikodym
derivative of $D b$ with respect to $|D b|$. Here the measure $|D b|$ on
$[0,T]\times\R^n$ is defined as
  $$\int \varphi(t,x)\,\d|Db|(t,x)=\int_0^T\!\!\int_{\R^n}\varphi(t,x)\,\d|D b_t|(x)\,\d t,
  \quad \mbox{for all } \varphi\in C_c([0,T]\times\R^n).$$
$Db,\, |D^a b|$ and $|D^s b|$ are defined similarly. We remark that Ambrosio
defined in \cite{Ambrosio04} the matrix $\tilde M_t$ to be the Radon-Nikodym
derivative of $D^s b$ with respect to $|D^s b|$; however, when restricted
on the support of $|D^s b|$, we have $M_t(x)=\tilde M_t(x)$ for $|D^s b|$-a.e.
$(t,x)\in[0,T]\times\R^n$.

\begin{lemma}[Commutator estimate II]\label{sect-2-lem.1}
For any compact set $K\subset(0,T)\times\R^n$, we have
  \begin{equation}\label{sect-2-lem.1.1}
  \varlimsup_{\ee\da0}\int_K |Q_{2,\ee}|\,\d x\d t
  \leq \|p\|_\infty\int_K\Lambda(M_t(x),\rho)\,\d|D^s b|(t,x)
  +\|p\|_\infty(n+I(\rho))|D^a b|(K)
  \end{equation}
and
  \begin{equation}\label{sect-2-lem.1.2}
  \varlimsup_{\ee\da0}\int_K |Q_{2,\ee}|\,\d x\d t
  \leq \|p\|_\infty I(\rho)|D^s b|(K),
  \end{equation}
where for a matrix $M$ and $\rho\in C^\infty_c(\R^n)$,
  $$\Lambda(M,\rho)=\int_{\R^n}|\<Mz,\nabla\rho(z)\>|\,\d z,\quad
  I(\rho)=\int_{\R^n}|z|\cdot|\nabla\rho(z)|\,\d z.$$
\end{lemma}

Now we turn to the error terms concerning the diffusion coefficient
$\sigma$. The arguments are similar to those in \cite[Subsection
5.3]{LeBrisLions08}. Applying again Lemma \ref{commutator-estimate},
we have
  \begin{equation}\label{R_epsilon}
  R_{k,\ee}=[\rho_\ee,\sigma_{jk}\partial_j](p)\stackrel{\ee\ra0}\lra 0
  \quad \mbox{in } L^2\big([0,T],L^2_{loc}(\R^n)\big),
  \end{equation}
due to the facts that $\sigma_{jk}\in
L^2\big([0,T],W^{1,2}_{loc}(\R^n)\big)$ and $p\in
L^\infty\big([0,T],L^\infty(\R^n)\big)$. Next, since both
$\partial_i\sigma_{ik}$ and $\sigma_{jk}\partial_j p$ belong to $
L^2\big([0,T],L^2_{loc}(\R^n) \big)$,
  \begin{equation}\label{S_epsilon}
  S_\ee=[\rho_\ee,\partial_i\sigma_{ik}](\sigma_{jk}\partial_j p)
  \stackrel{\ee\ra0}\lra 0  \quad \mbox{in } L^1\big([0,T],L^1_{loc}(\R^n)\big).
  \end{equation}
Finally, as $\sigma_{ik}\in L^2\big([0,T],W^{1,2}_{loc}(\R^n)\big),
\, i=1,\cdots,n$, we have
  \begin{equation}\label{T_epsilon}
  T_\ee=[\rho_\ee,\sigma_{ik}\partial_i](\sigma_{jk}\partial_j p)
  \stackrel{\ee\ra0}\lra 0  \quad \mbox{in } L^1\big([0,T],L^1_{loc}(\R^n)\big).
  \end{equation}

We denote by
  $$U_\ee=-Q_{1,\ee}+\frac12(S_\ee+T_\ee);$$
then the estimates \eqref{Q1}, \eqref{S_epsilon} and
\eqref{T_epsilon} lead to
  \begin{equation}\label{U_epsilon}
  U_\ee \stackrel{\ee\ra0}\lra 0 \quad \mbox{in }
  L^1\big([0,T],L^1_{loc}(\R^n)\big).
  \end{equation}
The equation \eqref{FPE_divergence_error} can be rewritten as
follows:
  \begin{equation}\label{FPE_divergence_error.1}
  \partial_t p_\ee+\partial_i(p_\ee b_i)-\frac12\partial_i(\sigma_{ik}\sigma_{jk}
  \partial_j p_\ee)=U_\ee-Q_{2,\ee}+\frac12\partial_i(\sigma_{ik}R_{k,\ee}).
  \end{equation}
Notice that $p_\ee$ is smooth with respect to the spatial variable.
By Remark \ref{sect-2-rem.1}, the commutator $R_{k,\ee}\in L^2\big([0,T],W^{1,2}_{loc}(\R^n)\big)$,
which together with the regularity assumptions on the coefficients $b$ and $\sigma$ tells
us that $\partial_t p_\ee\in L^1_{loc}\big([0,T]\times\R^n\big)$.
Therefore $p_\ee\in W^{1,1}_{loc}\big([0,T]\times\R^n\big)$ and we can apply
the standard chain rule in Sobolev spaces. For $\beta\in C^2(\R)$, one has
  \begin{align*}
  &\partial_t \beta(p_\ee)+\partial_i(\beta(p_\ee) b_i)-\frac12\partial_i(\sigma_{ik}\sigma_{jk}
  \partial_j \beta(p_\ee))\cr
  &\hskip6pt =\beta'(p_\ee)\partial_t p_\ee +b_i\beta'(p_\ee)\partial_i p_\ee
  +\beta(p_\ee) \div(b)-\frac12\partial_i\big(\sigma_{ik}\sigma_{jk}\beta'(p_\ee)\partial_j p_\ee\big)\cr
  &\hskip6pt =\beta'(p_\ee)\partial_t p_\ee +\beta'(p_\ee)\partial_i(p_\ee b_i)
  +\big(\beta(p_\ee)-p_\ee\beta'(p_\ee)\big)\div(b)-\frac12 \beta'(p_\ee)\partial_i(\sigma_{ik}\sigma_{jk}
  \partial_j p_\ee)\cr
  &\hskip16pt -\frac12(\sigma_{ik}\sigma_{jk} \partial_j p_\ee)\beta''(p_\ee)\partial_i p_\ee\cr
  &\hskip6pt =\beta'(p_\ee)\Big[\partial_t p_\ee+\partial_i(p_\ee b_i)-\frac12\partial_i(\sigma_{ik}\sigma_{jk}
  \partial_j p_\ee)\Big]+\big(\beta(p_\ee)-p_\ee\beta'(p_\ee)\big)\div(b)
  -\frac12\beta''(p_\ee)|\sigma^\ast \nabla p_\ee|^2.
  \end{align*}
By \eqref{FPE_divergence_error.1}, we obtain
  \begin{eqnarray}\label{FPE_divergence_error.2}
  &&\partial_t \beta(p_\ee)+\partial_i(\beta(p_\ee) b_i)-\frac12\partial_i(\sigma_{ik}\sigma_{jk}
  \partial_j \beta(p_\ee))-\big(\beta(p_\ee)-p_\ee\beta'(p_\ee)\big)\div(b)
  +\frac12\beta''(p_\ee)|\sigma^\ast \nabla p_\ee|^2\cr
  &&\hskip6pt = \beta'(p_\ee)\Big[U_\ee-Q_{2,\ee}+\frac12\partial_i(\sigma_{ik}R_{k,\ee})\Big].
  \end{eqnarray}

In order to prove the uniqueness of solutions to \eqref{FPE_divergence},
we shall use the technique of renormalized solutions, a notion which was introduced
by Di Perna and Lions in \cite[Section II.3]{DiPernaLions89}
(see also \cite[Definition 4.1]{Ambrosio08} and \cite[Definition 4.9]{Figalli}).

\begin{definition}[Renormalized solution]\label{sect-2-def}
Let $\sigma:[0,T]\times\R^n\ra \R^m\otimes\R^n$ and $b:[0,T]\times\R^n\ra \R^n$ be such that
(i) $\sigma\in L^2_{loc}([0,T]\times\R^n)$,
(ii) $b,\,\div(b)\in L^1_{loc}([0,T]\times\R^n)$.
We say that a solution $p$ to \eqref{FPE_divergence} is a renormalized solution if for any
$\beta\in C^2(\R)$, the following equation holds in the distributional sense:
  \begin{equation}\label{sect-2-def.1}
  \partial_t \beta(p)+\partial_i(\beta(p) b_i)-\frac12\partial_i(\sigma_{ik}\sigma_{jk}
  \partial_j \beta(p))-\big(\beta(p)-p\beta'(p)\big)\div(b)
  +\frac12\beta''(p)|\sigma^\ast \nabla p|^2=0.
  \end{equation}
\end{definition}

Now we show that any weak solution in the space $X_2$ of \eqref{FPE_divergence}
is renormalizable, provided the conditions of Theorem \ref{Main_thm} are satisfied.
The main idea of the arguments are similar to the proof of \cite[Theorem 3.5]{Ambrosio04}
(see in particular Step 3 therein).

\begin{theorem}[Renormalization property]\label{sect-2-thm}
Under the conditions of Theorem \ref{Main_thm}, any weak solution
  $$p\in X_2=\big\{p\in L^\infty\big([0,T],L^2\cap L^\infty\big);
  \sigma^\ast\nabla p\in L^2\big([0,T],L^2\big)\big\}$$
of \eqref{FPE_divergence} is also a renormalized solution.
\end{theorem}

\noindent{\bf Proof.} We have to show that as $\ee\da0$, all the terms
on the left hand side of the equation \eqref{FPE_divergence_error.2}
converge in the distributional sense to the corresponding ones in \eqref{sect-2-def.1},
while the limit of the right hand side is 0.
We split the proof into two steps: in the first step we show the convergences
of all the terms except the one involving $Q_{2,\ee}$, while in the second step
we focus on the term $\beta'(p_\ee)Q_{2,\ee}$.

{\bf Step 1.}
Since $p$ is essentially bounded and $\beta\in C^2(\R)$, $\beta$ and its derivatives are uniformly
continuous on the interval $I_p=[-\|p\|_\infty,\|p\|_\infty]$.
Notice also that $\|p_\ee\|_\infty\leq \|p\|_\infty$. As the proofs
are similar, we only illustrate the convergences of
  $$\partial_i(\sigma_{ik}\sigma_{jk}\partial_j \beta(p_\ee))\ra
  \partial_i(\sigma_{ik}\sigma_{jk}\partial_j \beta(p))\quad
  \mbox{and} \quad \beta'(p_\ee)\partial_i(\sigma_{ik}R_{k,\ee})\ra 0$$
as $\ee\da0$ in the distributional sense.

First for any $\varphi\in C_c^\infty\big([0,T)\times\R^n\big)$, we have
by the integration by parts formula that
  \begin{align*}
  &\bigg|\int_0^T\!\!\int_{\R^n}\varphi\, \partial_i(\sigma_{ik}\sigma_{jk}\partial_j \beta(p_\ee))\,\d x\d t
  -\int_0^T\!\!\int_{\R^n}\varphi\, \partial_i(\sigma_{ik}\sigma_{jk}\partial_j \beta(p))\,\d x\d t\bigg|\cr
  &\hskip6pt \leq \int_0^T\!\!\int_{\R^n}\big|\beta'(p_\ee)\<\sigma^\ast\nabla\varphi,\sigma^\ast\nabla p_\ee\>
  -\beta'(p)\<\sigma^\ast\nabla\varphi,\sigma^\ast\nabla p\>\big|\,\d x\d t.
  \end{align*}
The above quantity is dominated by the sum of
  $$J_{1,\ee}:=\|\beta'(p_\ee)\|_{\infty}\|\nabla\varphi\|_\infty\int_K
  |\sigma|\cdot|\sigma^\ast\nabla p_\ee-\sigma^\ast\nabla p|\,\d x\d t$$
and
  $$J_{2,\ee}:=\|\nabla\varphi\|_\infty\int_K |\beta'(p_\ee)-\beta'(p)|
  \cdot|\sigma|\cdot|\sigma^\ast\nabla p|\,\d x\d t,$$
where the compact set $K:=\supp(\varphi)\subset [0,T]\times \R^n$ is the support of $\varphi$.
By Cauchy's inequality, we have
  \begin{align*}
  J_{1,\ee}&\leq C\bigg( \int_K |\sigma|^2\d x\d t \bigg)^{\frac12}
  \bigg( \int_K |\sigma^\ast\nabla p_\ee-\sigma^\ast\nabla p|^2\d x\d t \bigg)^{\frac12}
  \ra0
  \end{align*}
as $\ee\da0$. In this paper $C$ denotes the constants whose values have no
importance and may change from line to line. Next since $p\in L^\infty([0,T],L^2(\R^n))
\subset L^2([0,T],L^2(\R^n))$,
$p_\ee$ tends to $p$ in the latter space as $\ee\da0$; consequently,
$\beta'(p_\ee)$ converges to $\beta'(p)$ in measure. Moreover,
  $$|\beta'(p_\ee)-\beta'(p)|\cdot|\sigma|\cdot|\sigma^\ast\nabla p|
  \leq \big(\|\beta'(p_\ee)\|_\infty+\|\beta'(p)\|_\infty\big)|\sigma|\cdot|\sigma^\ast\nabla p|
  \leq C|\sigma|\cdot|\sigma^\ast\nabla p|$$
which is integrable on $K$. Therefore the dominated convergence theorem tells us that
  $$\lim_{\ee\da0}J_{2,\ee}=0.$$
To sum up, we conclude that when $\ee\da0$, $\partial_i(\sigma_{ik}\sigma_{jk}\partial_j \beta(p_\ee))$
converges to $\partial_i(\sigma_{ik}\sigma_{jk}\partial_j \beta(p))$ in the distributional sense.

Now we consider the limit $\beta'(p_\ee)\partial_i(\sigma_{ik}R_{k,\ee})\ra 0$.
For $\varphi\in C_c^\infty\big([0,T)\times\R^n\big)$, again by integrating by parts, one has
  \begin{align*}
  &\int_0^T\!\!\int_{\R^n}\varphi\, \beta'(p_\ee)\partial_i(\sigma_{ik}R_{k,\ee})\,\d x\d t\cr
  &\hskip6pt =-\int_0^T\!\!\int_{\R^n}\beta'(p_\ee)\<\sigma^\ast\nabla\varphi, R_\ee\>\,\d x\d t
  -\int_0^T\!\!\int_{\R^n}\varphi\, \beta''(p_\ee)\<\sigma^\ast\nabla p_\ee, R_\ee\>\,\d x\d t\cr
  &\hskip6pt =:-J_{3,\ee}-J_{4,\ee},
  \end{align*}
where $R_\ee=(R_{1,\ee},\cdots,R_{m,\ee})$. Note that
  \begin{align*}
  |J_{3,\ee}|&\leq \|\beta'(p_\ee)\|_\infty\|\nabla\varphi\|_\infty
  \int_K |\sigma|\cdot|R_\ee|\,\d x\d t\cr
  &\leq C\bigg( \int_K |\sigma|^2\,\d x\d t \bigg)^{\frac12}
  \bigg( \int_K |R_\ee|^2\,\d x\d t \bigg)^{\frac12}.
  \end{align*}
Therefore by \eqref{R_epsilon}, we get $\lim_{\ee\da0}J_{3,\ee}=0$.
Since $\beta''(p_\ee)$ and $\|\sigma^\ast\nabla p_\ee\|_{L^2([0,T]\times\R^n)}$
are bounded, uniformly in $\ee>0$,
we can show that $\lim_{\ee\da0}J_{4,\ee}=0$ in the same way.

Summing up these discussions, we conclude that as $\ee\da0$,
all the terms in \eqref{FPE_divergence_error.2},
except $-\beta'(p_\ee)Q_{2,\ee}$, converge in the distributional sense.
We define the ``defect'' measure
  \begin{align}\label{sect-2-thm.1}
  \mu&:=\partial_t \beta(p)+\partial_i(\beta(p) b_i)-\frac12\partial_i(\sigma_{ik}\sigma_{jk}
  \partial_j \beta(p))-\big(\beta(p)-p\beta'(p)\big)\div(b)+\frac12\beta''(p)|\sigma^\ast \nabla p|^2
  \end{align}
(notice that $\beta'(p_\ee)\big[U_\ee+\frac12\partial_i(\sigma_{ik}R_{k,\ee})\big]\ra0$
as $\ee$ tends to 0).

{\bf Step 2.} Now we deal with the term $-\beta'(p_\ee)Q_{2,\ee}$.
Let $Q_\rho$ be one of the weak limit points of $|\beta'(p_\ee)Q_{2,\ee}|$ in the sense of measure
(by Lemma \ref{sect-2-lem.1}, such an accumulating point exists since $|\beta'(p_\ee)Q_{2,\ee}|$
is bounded in $L^1_{loc}$). The measure $Q_\rho$ depends on the convolution kernel $\rho$,
but the ``defect'' measure $\mu$ defined in \eqref{sect-2-thm.1}
is independent of $\rho$ and satisfies $|\mu| \leq Q_\rho$.
Thus we get
  \begin{equation}\label{sect-2-thm.2}
  |\mu| \leq Q_\rho\quad \mbox{for all }\rho\in\K.
  \end{equation}
We deduce from \eqref{sect-2-lem.1.2} that $Q_\rho$ is absolutely continuous
with respect to $|D^s b|$, which together
with  \eqref{sect-2-lem.1.1} gives us
  $$Q_\rho\leq \|p\|_\infty \Lambda(M_\cdot(\cdot),\rho)|D^s b|.$$
Thus by \eqref{sect-2-thm.2} we obtain
  \begin{equation*}
  |\mu| \leq \|p\|_\infty \Lambda(M_\cdot(\cdot),\rho)|D^s b|,\quad \mbox{for all }\rho\in\K.
  \end{equation*}
Denote by $g$ the Radon-Nikodym density of $|\mu|$ with respect to $|D^s b|$. Then for each
$\rho\in\K$,
  $$g(t,x)\leq \|p\|_\infty \Lambda(M_t(x),\rho)\quad \mbox{for } |D^s b|\mbox{-a.e. } (t,x).$$
Let $\D$ be a countable dense subset of $\K$ with respect to the norm $W^{1,1}(B(1))$.
We have
  \begin{equation}\label{sect-2-thm.4}
  g(t,x)\leq \|p\|_\infty \inf_{\rho\in\D}\Lambda(M_t(x),\rho)
  \quad \mbox{for } |D^s b|\mbox{-a.e. } (t,x).
  \end{equation}
For $|D^s b|$-a.e. $(t,x)$ fixed, we deduce from the definition of $\Lambda(M_t(x),\rho)$
that the mapping $\K\ni \rho\ra\Lambda(M_t(x),\rho)$ is continuous with respect to
the $W^{1,1}(B(1))$ norm. Therefore
  $$\inf_{\rho\in\D}\Lambda(M_t(x),\rho)=\inf_{\rho\in\K}\Lambda(M_t(x),\rho).$$
Now by Alberti's rank one structure of $M_t(x)$ (cf. \cite[Theorem 2.3]{Ambrosio04}), we conclude that
the above infimum is 0 for $|D^s b|$-a.e. $(t,x)\in[0,T]\times\R^n$ (see \cite[Lemma 3.3]{Ambrosio04}).
Therefore $g=0$ $|D^s b|$-a.e. by \eqref{sect-2-thm.4}. As a result, the measure $\mu=0$.
By \eqref{sect-2-thm.1}, we arrive at
  \begin{equation}\label{sect-2-thm.5}
  \partial_t \beta(p)+\partial_i(\beta(p) b_i)-\frac12\partial_i(\sigma_{ik}\sigma_{jk}
  \partial_j \beta(p))-\big(\beta(p)-p\beta'(p)\big)\div(b)
  +\frac12\beta''(p)|\sigma^\ast \nabla p|^2=0.
  \end{equation}
The proof is complete.  \fin

\medskip

Finally we are ready to prove the uniqueness of solutions to \eqref{FPE.1}.

\medskip

{\noindent \bf Proof of Theorem \ref{Main_thm}.} Since the equation \eqref{FPE.1}
is linear, it is enough to show that if $p$ is a weak solution such that $p(0)\equiv 0$,
then $p(t)\equiv 0$ for all $t\geq0$. By Theorem \ref{sect-2-thm}, we
know that for any $\beta\in C^2(\R)$, \eqref{sect-2-thm.5} holds in the
distributional sense. If we choose $\beta(s)=s^2$ for $s\in\R$, then it becomes
  \begin{equation}\label{proof.1}
  \partial_t (p^2)+\partial_i(p^2 b_i)-\partial_i(\sigma_{ik}\sigma_{jk}
  p\,\partial_j p)+p^2\div(b)\leq 0.
  \end{equation}

We take $\phi\in C_c^\infty(\R^n,[0,1])$ satisfying
$\phi|_{B(1)}\equiv 1$ and $\supp(\phi)\subset B(2)$. Define a
nonnegative smooth cut-off function $\phi_R=\phi(\frac\cdot R)$ for
$R>0$; then $\nabla\phi_R= \frac1R\nabla\phi(\frac\cdot R)$. To
simplify the notations, we will write $\int_{\R^n}f$ for the
integral of the function $f$ on $\R^n$ with respect to the Lebesgue
measure. Multiplying the inequality \eqref{proof.1}
with $\phi_R$ and integrating by parts on $\R^n$, we obtain
  \begin{align}\label{proof.2}
  \frac\d{\d t}\int_{\R^n}p^2\phi_R
  \leq\int_{\R^n}p^2\<b,\nabla\phi_R\>-\int_{\R^n}p\,\<\sigma^\ast\nabla\phi_R,\sigma^\ast\nabla p\>
  -\int_{\R^n} p^2\phi_R\div(b).
  \end{align}

Now we estimate the three terms on the right hand side. First
  \begin{align*}
  \bigg|\int_{\R^n}p^2\< b,\nabla\phi_R\>\bigg|
  &=\bigg|\int_{\{R\leq|x|\leq 2R\}}p^2\big\< b,\frac1R \nabla\phi(\frac\cdot
  R)\big\>\bigg|\cr
  &\leq C\|\nabla\phi\|_{L^\infty}
  \int_{\{|x|\geq R\}}p^2\frac{|b|}{1+|x|}.
  \end{align*}
Since $\frac
b{1+|x|}\in \big(L^1([0,T],L^1+L^\infty(\R^n))\big)^n$, there are
two vector fields $b_1,b_2$ such that $b=b_1+b_2$ and $\frac{|b_1|}{1+|x|}\in
L^1([0,T],L^1(\R^n))$, $\frac{|b_2|}{1+|x|}\in
L^1([0,T],L^\infty(\R^n))$. Then we have
  \begin{align}\label{proof.3}
  \bigg|\int_{\R^n}p^2\< b,\nabla\phi_R\>\bigg|
  &\leq C\|p\|^2_{L^\infty([0,T],L^\infty)}\int_{\{|x|\geq R\}}\frac{|b_1(t)|}{1+|x|}
  +C\bigg\|\frac{b_2(t)}{1+|x|}\bigg\|_{L^\infty}\int_{\{|x|\geq R\}} p^2.
  \end{align}
Note that as a function of $t\in[0,T]$, the right hand side of
\eqref{proof.3} is dominated by
  $$C\bigg\|\frac{b_1(t)}{1+|x|}\bigg\|_{L^1}
  +C\|p\|^2_{L^\infty([0,T],L^2)}\bigg\|\frac{b_2(t)}{1+|x|}\bigg\|_{L^\infty},$$
and the latter is an integrable function of $t\in[0,T]$. Furthermore, for
a.e. $t\in[0,T]$, the right hand side of  \eqref{proof.3} tends to 0
as $R\ra\infty$. Therefore by the dominated convergence theorem,
  \begin{equation}\label{proof.4}
  \lim_{R\ra\infty}\int_0^T\!\!\int_{\R^n}p^2\<b,\nabla\phi_R\>=0。
  \end{equation}

Next, in a similar way,
  \begin{align*}
  \bigg|\int_{\R^n}p\,\< \sigma^\ast\nabla
  p,\sigma^\ast\nabla\phi_R\>\bigg|
  &\leq \int_{\{R\leq|x|\leq 2R\}}|p |\cdot |\sigma^\ast\nabla
  p|\cdot\frac1R \|\nabla\phi\|_{L^\infty}|\sigma|\cr
  &\leq C\|\nabla\phi\|_{L^\infty}
  \int_{\{|x|\geq R\}}|p|\cdot|\sigma^\ast\nabla p|\frac{|\sigma|}{1+|x|}.
  \end{align*}
By the assumptions on $\sigma$, we can split it into two
matrix-valued functions $\sigma_1$ and $\sigma_2$, such that
  \begin{equation*}
  \frac{\sigma_1}{1+|x|}\in L^2([0,T],L^2(\R^n))\quad \mbox{and}
  \quad \frac{\sigma_2}{1+|x|}\in L^2([0,T],L^\infty(\R^n)).
  \end{equation*}
By Cauchy's inequality,
  \begin{align}\label{proof.5}
  \bigg|\int_{\R^n}p\,\< \sigma^\ast\nabla
  p,\sigma^\ast\nabla\phi_R\>\bigg|
  &\leq C\|p\|_{L^\infty}\int_{\{|x|\geq R\}}|\sigma^\ast\nabla p|\frac{|\sigma_1|}{1+|x|}
  +C\int_{\{|x|\geq R\}}|p|\cdot |\sigma^\ast\nabla
  p|\frac{|\sigma_2|}{1+|x|}\cr
  &\leq C\|\sigma^\ast\nabla p\|_{L^2}\bigg(\int_{\{|x|\geq
  R\}}\frac{|\sigma_1(t)|^2}{(1+|x|)^2}\bigg)^{1/2}\cr
  &\hskip11pt +C\bigg\|\frac{\sigma_2(t)}{1+|x|}\bigg\|_{L^\infty}\int_{\{|x|\geq R\}}|p|\cdot |\sigma^\ast\nabla
  p|.
  \end{align}
It is easy to see that the right hand side of \eqref{proof.5} is
dominated by
  $$C\|\sigma^\ast\nabla p\|_{L^2}\bigg\|\frac{\sigma_1(t)}{1+|x|}\bigg\|_{L^2}
  +C\|p\|_{L^\infty([0,T],L^2)}\|\sigma^\ast\nabla
  p\|_{L^2}\bigg\|\frac{\sigma_2(t)}{1+|x|}\bigg\|_{L^\infty},$$
which, by Cauchy's inequality and the properties of $\sigma_1,\sigma_2$ and $p$, is
integrable with respect to $t\in[0,T]$. Moreover, for a.e.
$t\in[0,T]$, the right hand side of \eqref{proof.5} vanishes as
$R\ra\infty$. Consequently, the dominated convergence theorem leads to
  \begin{equation}\label{proof.6}
  \lim_{R\ra\infty}\int_0^T\!\!\int_{\R^n}p\,\< \sigma^\ast\nabla
  p,\sigma^\ast\nabla\phi_R\>=0.
  \end{equation}

Finally it is clear that
  \begin{equation}\label{proof.7}
  -\int_{\R^n}p^2\phi_R\,\div(b)\leq \int_{\R^n}p^2\phi_R\, [\div(b)]^-
  \leq \|[\div(b)]^-\|_{L^\infty}\int_{\R^n}p^2.
  \end{equation}
Now we integrate the equation \eqref{proof.2} in time from $0$
to $t$ (note that $p(0)\equiv 0$) and obtain
  $$\int_{\R^n}p^2\phi_R\leq\int_0^T\!\!\int_{\R^n}p^2\<b,\nabla\phi_R\>
  -\int_0^T\!\!\int_{\R^n}p\,\<\sigma^\ast\nabla\phi_R,\sigma^\ast\nabla p\>
  -\int_0^T\!\!\int_{\R^n} p^2\phi_R\div(b).$$
By \eqref{proof.4} and \eqref{proof.6}, for any $\eta>0$ there exists
$R_0>0$ such that for any $R>R_0$, we have
  $$\bigg|\int_0^T\!\!\int_{\R^n}p^2\<b,\nabla\phi_R\>
  -\int_0^T\!\!\int_{\R^n}p\,\<\sigma^\ast\nabla\phi_R,\sigma^\ast\nabla p\>\bigg|\leq \eta.$$
Taking into consideration the inequality \eqref{proof.7}, we get
  $$\int_{\R^n}p^2\phi_R\leq \eta+\int_0^T\|[\div(b)]^-\|_{L^\infty}\int_{\R^n}p^2.$$
First letting $R\ua+\infty$ and then $\eta\da0$, we finally get
  $$\int_{\R^n}p^2 \leq \int_0^T\|[\div(b)]^-\|_{L^\infty}\int_{\R^n}p^2.$$
Since $[\div(b)]^-\in L^1\big([0,T],L^\infty(\R^n)\big)$, we have $\int_{\R^n}p^2=0$.
This shows $p=0$ a.e. and the uniqueness is proved. \fin

\begin{remark}\label{sect-2-rem.2} \rm
Now we briefly discuss the reason why we are unable to deal with the
diffusion coefficients $\sigma$ of BV regularity. Recall the
definition of $S_\ee$ in \eqref{S-epsilon_T-epsilon}:
  $$S_\ee=[\rho_\ee,\partial_i\sigma_{ik}](\sigma_{jk}\partial_jp)
  =\rho_\ee\ast\big((\sigma_{jk}\partial_jp)\partial_i\sigma_{ik}\big)
  -\big(\rho_\ee\ast(\sigma_{jk}\partial_jp)\big)\partial_i\sigma_{ik}.$$
For simplicity, we assume that the functions are time-independent
and denote by $g=\sigma_{jk}\partial_jp\in L^2(\R^n)$.
Then
  \begin{equation*}
  S_\ee=\rho_\ee\ast(g\partial_i\sigma_{ik})
  -(\rho_\ee\ast g)\partial_i\sigma_{ik}.
  \end{equation*}
If $\sigma_{ik}$ has only $\BV_{loc}$ regularity, then
$\mu:=\partial_i\sigma_{ik}$ is a locally finite Radon measure on $\R^n$.
We have the decomposition $\mu=D^a \mu+D^s \mu$, where $D^a\mu\ll\d x$
and $D^s\mu \perp\d x$. In the case $D^s\mu\neq0$, since $g\in L^2(\R^n)$
is not continuous, the product $g D^s\mu$
is sensitive to the modification of $g$ in Lebesgue negligible sets.
Therefore $S_\ee$ is not a well-defined object on $\R^n$. \fin
\end{remark}

\end{document}